\documentclass[11pt]{article}
\usepackage{mathrsfs}
\usepackage{amsfonts,amsmath,dsfont,amssymb,amsthm,stmaryrd}
\usepackage[fit]{truncate}
\usepackage{amsopn, graphicx}
\usepackage{hyperref}
\hypersetup{colorlinks=true,pdftitle="",pdftex}

\newcommand{\bc}{\begin{center}}
\newcommand{\ec}{\end{center}}
\newcommand{\bt}{\begin{tabular}}
\newcommand{\et}{\end{tabular}} 
\newcommand{\bea}{\begin{eqnarray}}
\newcommand{\eea}{\end{eqnarray}}
\newcommand{\bean}{\begin{eqnarray*}}
\newcommand{\eean}{\end{eqnarray*}}

\newcommand{\ba}{\begin{array}}
\newcommand{\ea}{\end{array}}

\def\be{\begin{eqnarray}}
\def\ee{\end{eqnarray}}
\def\ben{\begin{eqnarray*}}
\def\een{\end{eqnarray*}}

\newcommand{\ra} {\rightarrow}

\newcommand{\RL}{{\mathbb R}}

\newcommand{\lam}{\lambda}

\def\elabel#1{\label{e:#1}}

\def\sq{$\Box$}

\def\qed{\ifmmode\sq\else{\unskip\nobreak\hfil
\penalty50\hskip1em\null\nobreak\hfil\sq
\parfillskip=0pt\finalhyphendemerits=0\endgraf}\fi\par\medbreak}

\newsavebox{\junk}
\savebox{\junk}[1.6mm]{\hbox{$|\!|\!|$}}

\def\til={{\widetilde =}}

 \def\eq#1/{(\ref{#1})}

\def\eq#1/{(\ref{e:#1})}

\newcommand{\beqn}[1]{\notes{#1}%
\begin{eqnarray} \elabel{#1}}

\newcommand{\eeqn}{\end{eqnarray} }

\newcommand{\beq}[1]{\notes{#1}%
\begin{equation}\elabel{#1}}

\newcommand{\eeq}{\end{equation}} 

\def\bdes{\begin{description}}
\def\edes{\end{description}}

\def\notes#1{}

\newtheorem{conj}{Conjecture}[section]
\newtheorem{thm}{Theorem}[section]

\newcommand{\conv}{\text{conv}}

\newcommand{\R}{\RL}

\def\R{{\bf R}}

\def\bee{\begin{eqnarray*}}
\def\ene{\end{eqnarray*}}

\newcommand{\vol}{\mathrm{Vol}}

\newcommand{\N}{\mathbb{N}}

\begin{document}
\title{Do Minkowski averages get progressively more convex?}
\author{Matthieu Fradelizi\thanks{Laboratoire d'Analyse et de Math\'ematiques Appliqu\'ees UMR 8050,
Universit\'e Paris-Est Marne-la-Vall\'ee, 5 Bd Descartes, Champs-sur-Marne,
77454 Marne-la-Vall\'ee Cedex 2, France. Supported in part by the Agence Nationale de la Recherche, project GeMeCoD (ANR 2011 BS01 007 01).
Email: {\tt matthieu.fradelizi@u-pem.fr}}, 
Mokshay Madiman\thanks{Department of Mathematical Sciences, University of Delaware, 501 Ewing Hall,
Newark, DE 19716, USA. Supported in part by the U.S. National Science Foundation through grants DMS-1409504 (CAREER) and CCF-1346564.
Email: {\tt madiman@udel.edu}}, 
Arnaud Marsiglietti\thanks{Institute for Mathematics and its Applications, University of Minnesota,
207 Church Street SE, 434 Lind Hall, Minneapolis, MN 55455, USA. Supported in part by the Institute for Mathematics and its Applications with funds provided by the National Science Foundation. Email: {\tt arnaud.marsiglietti@ima.umn.edu}} 
and Artem Zvavitch\thanks{Department of Mathematical Sciences, Kent State University, Kent, OH 44242, USA. 
Supported in part by the U.S. National Science Foundation Grant DMS-1101636 and the Simons Foundation. Email: {\tt zvavitch@math.kent.edu}}}
\date{}
\maketitle


\begin{abstract}
Let us define, for a compact set $A \subset \R^n$, the Minkowski  averages of $A$:
$$ A(k) = \left\{\frac{a_1+\cdots +a_k}{k} : a_1, \ldots, a_k\in A\right\}=\frac{1}{k}\Big(\underset{k\ {\rm times}}{\underbrace{A + \cdots + A}}\Big). $$
We study the monotonicity of the convergence of $A(k)$ towards the convex hull of $A$,  when 
considering the Hausdorff distance, the volume deficit and a non-convexity index of Schneider as measures of convergence.
For the volume deficit, we show that monotonicity fails in general, thus disproving a conjecture of Bobkov, Madiman and Wang.
For Schneider's non-convexity index, we prove that a strong form of monotonicity holds, and for the Hausdorff distance, we establish 
that the sequence is eventually nonincreasing.
\end{abstract}

\section{Introduction}
\label{sec:intro}

This note announces and proves some of the results obtained in \cite{FMMZ15}. Let us denote for a compact set $A \subset \R^n$ and for a positive integer $k$,
\begin{eqnarray}\label{defAk}
A(k) = \left\{\frac{a_1+\cdots +a_k}{k} : a_1, \ldots, a_k\in A\right\}=\frac{1}{k}\Big(\underset{k\ {\rm times}}{\underbrace{A + \cdots + A}}\Big).
\end{eqnarray}
Denoting by $\conv(A)$ the convex hull of $A$, and by 
$$
d(A):= \inf\{r>0: \conv(A)\subset A+rB_2^n\}
$$ 
the Hausdorff distance between a set $A$ and its convex hull, it is a classical fact (proved independently by \cite{Sta69,EG69} in 1969,
and often called the Shapley-Folkmann-Starr theorem) that 
$A(k)$ converges in Hausdorff distance to $\conv(A)$ as $k \to \infty$.
Furthermore \cite{Sta69,EG69} also determined the rate of convergence: it turns out that $d(A(k))=O(1/k)$ for any compact set $A$.
For sets of nonempty interior, this convergence of Minkowski averages to the convex hull can also be expressed in terms of the volume deficit $\Delta(A)$ of a compact set $A$ in $\R^n$, which is defined as:
$$ \Delta(A):= \vol_n(\conv(A)\setminus A) = \vol_n(\conv(A))-\vol_n(A), $$
where $\vol_n$ denotes Lebesgue measure in $\R^n$. It was shown by \cite{EG69} that if $A$ is compact with nonempty interior, then the volume deficit of $A(k)$ also converges to 0; more precisely, $\Delta(A(k))=O(1/k)$ for any compact set $A$ with nonempty interior.

Our original motivation came from a conjecture made by Bobkov, Madiman and Wang \cite{BMW11}:

\vspace{.1in}
\begin{conj}[\cite{BMW11}]\label{weakconj}
Let $A$ be a compact set in $\R^n$ for some $n\in\N$, and let $A(k)$ be defined as in \eqref{defAk}. Then  the sequence $\Delta(A(k))$ is non-increasing in $k$, or equivalently, $\{\vol_n(A(k))\}_{k \geq 1}$ is non-decreasing.
\end{conj}
\vspace{.1in}

We show that Conjecture \ref{weakconj} fails to hold in general, even for moderately high dimension.

\vspace{.1in}
\begin{thm}\label{voldef}
Conjecture \ref{weakconj} is false in $\R^n$ for $n\ge12$, and true for $\RL^1$.
\end{thm}
\vspace{.1in}
Notice that Conjecture \ref{weakconj} remains open for $1<n<12$. In particular, the arguments presented in this note do not seem to work.
In analogy with Conjecture \ref{weakconj}, we also consider whether one can have monotonicity of $\{c(A(k))\}_{k \geq 1}$, where $c$ is a non-convexity index defined by Schneider \cite{Sch75} as follows: 
$$ c(A) := \inf \{ \lambda\geq 0: A+\lambda\, \conv(A) \text{ is convex}  \}. $$
A nice property of Schneider's index is that it is affine-invariant, i.e., 
$c(TA+x)=c(A)$ for any nonsingular linear map $T$ on $\R^n$ and any $x\in\R^n$. 

Contrary to the volume deficit, we prove that Schneider's non-convexity index $c$ satisfies a strong kind of monotonicity in any dimension.

\vspace{.1in}
\begin{thm}\label{thm:c-quant-mono}
Let $A$ be a compact set in $\R^n$ and $k\in\N^*$. Then 
$$ c\left(A(k+1)\right)\le\frac{k}{k+1}c\left(A(k)\right). $$
\end{thm}
\vspace{.1in}

Finally, we also prove that eventually, for $k\ge c(A)$, the Hausdorff distance between $A(k)$ and $\conv(A)$ is also strongly decreasing.

\vspace{.1in}
\begin{thm}\label{thm:d-quant-mono}
Let $A$ be a compact set in $\R^n$ and $k\ge c(A)$ be an integer. Then 
$$ d\left(A(k+1)\right)\le\frac{k}{k+1}d\left(A(k)\right). $$
\end{thm}
\vspace{.1in}

Moreover, Schneider proved in \cite{Sch75} that $c(A)\le n$ for every compact subset $A$ of $\R^n$. It follows that the eventual monotonicity of the sequence $d\left(A(k)\right)$ holds true for $k\ge n$.

It is natural to ask what the relationship is in general between convergence of $c$, $\Delta$ and $d$ to 0,
for arbitrary sequences $(C_k)$ of compact sets. In fact, none of these 3 notions of approach to convexity
are comparable with each other in general.
To see why, observe that while $c$ is scaling-invariant, neither $\Delta$ nor $d$ are; so it is easy to
construct examples of sequences $(C_k)$ such that $c(C_k)\ra 0$ but $\Delta(C_k)$ and $d(C_k)$
remain bounded away from 0. The same argument enables us to construct examples of sequences $(C_k)$ such that $c(C_k)$ remain bounded away from 0, whereas $\Delta(C_k)$ and $d(C_k)$ converge to 0. Furthermore, $\Delta(C_k)$ remains bounded away from 0 for any sequence $C_k$ of finite sets, whereas $c(C_k)$ and $d(C_k)$ could converge to 0 if the finite sets form a finer and finer grid filling out a convex set.
An example where $\Delta(C_k)\ra 0$ but both $c(C_k)$ and $d(C_k)$ are bounded away from 0
is given by taking a 3-point set with 2 of the points getting arbitrarily closer but staying away from the third.
One can obtain further relationships between these measures of non-convexity if further conditions are imposed on the sequence $C_k$; details may be found in \cite{FMMZ15}.

The rest of this note is devoted to the examination of whether $A(k)$ becomes progressively more convex as $k$ increases,
when measured through the functionals $\Delta, d$ and $c$. The concluding section contains some additional discussion.

\section{The behavior of volume deficit}\label{sec:Delta}

We prove Theorem \ref{voldef} in this section. 
We start by constructing a counterexample to the conjecture in $\R^{n}$, for $n\geq 12$.
Let $F$ be a $p$-dimensional subspace of $\R^n$, where $p \in \{1,\ldots, n-1\}$. Let us consider $A=I_1 \cup I_2$, where $I_1 \subset F$ and $I_2 \subset F^{\perp}$
are convex sets, and $F^{\perp}$ denotes the orthogonal complement of $F$. One has
$$ A+A = 2I_1 \cup (I_1 \times I_2) \cup 2 I_2, $$
$$ A+A+A = 3I_1 \cup (2I_1 \times I_2) \cup (I_1 \times 2I_2) \cup 3 I_2. $$
Notice that
$$ \vol_n(A+A) = \vol_p(I_1)\vol_{n-p}(I_2), $$
$$ \vol_n(A+A+A) = \vol_p(I_1)\vol_{n-p}(I_2) (2^p + 2^{n-p} - 1). $$
Thus, $\vol_n(A(3)) \geq \vol_n(A(2))$ if and only if
\begin{eqnarray}\label{calcul}
2^p + 2^{n-p} - 1 \geq \left( \frac{3}{2} \right)^n.
\end{eqnarray}
Notice that inequality~(\ref{calcul}) does not hold when $n \geq 12$ and $p=\lceil \frac{n}{2} \rceil$.

\begin{figure}[h!]\label{fig:dcr}

\begin{center}
\includegraphics[scale=0.45]{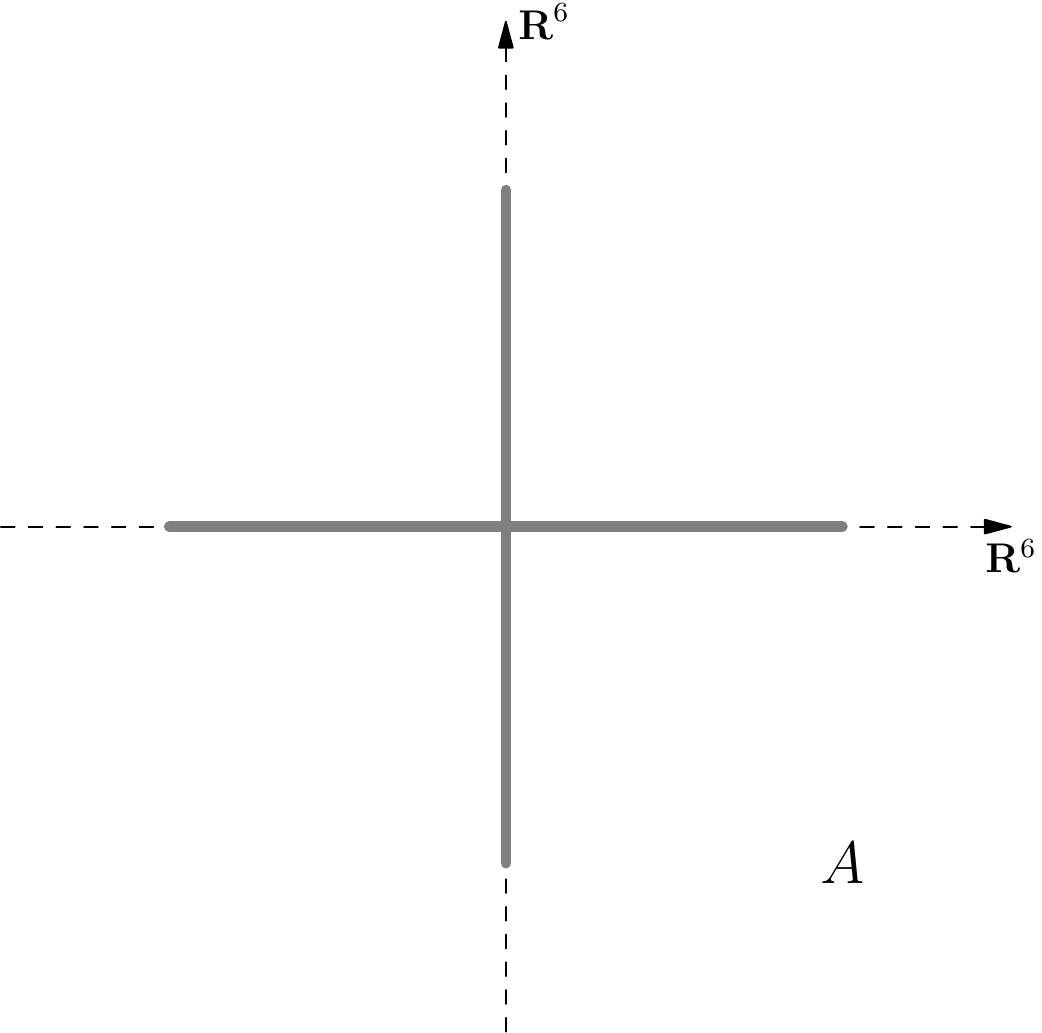} \,
\includegraphics[scale=0.45]{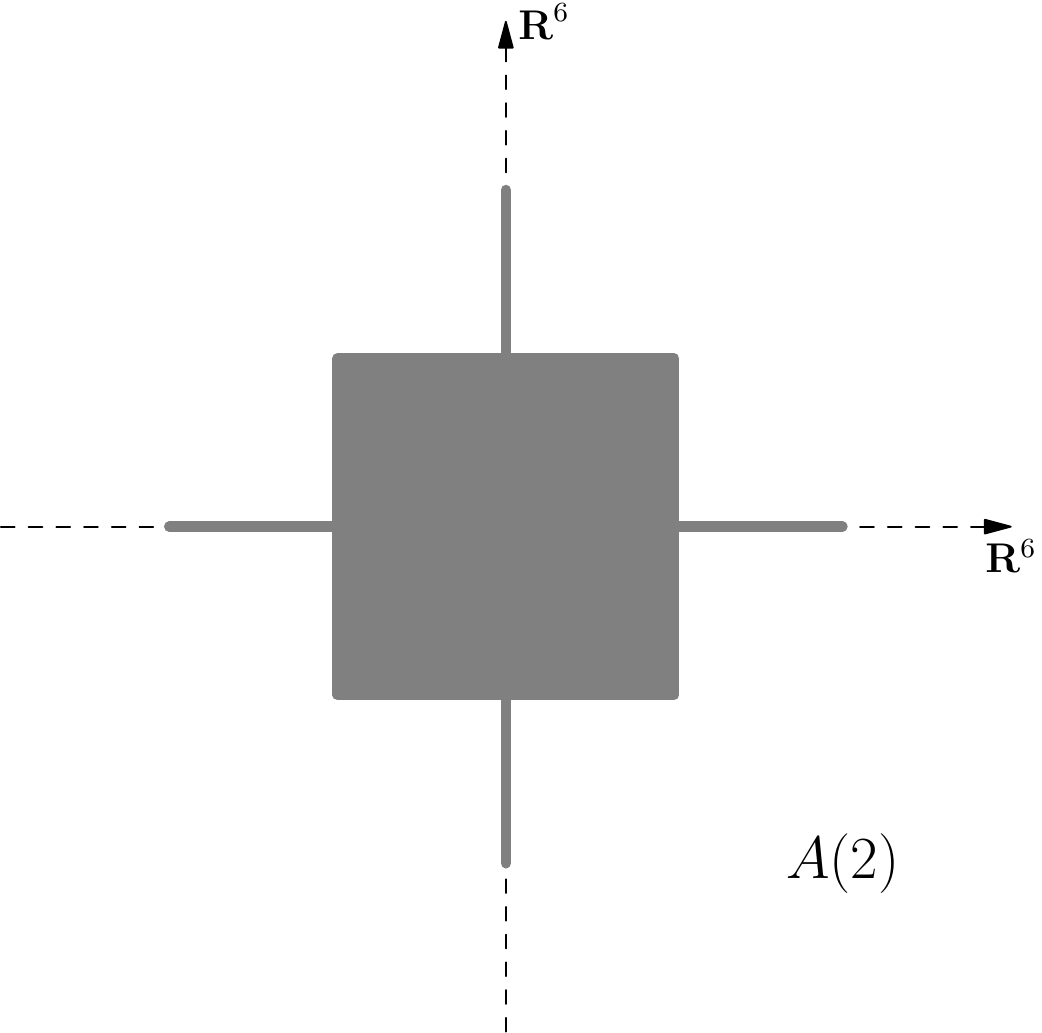}\,
\includegraphics[scale=0.45]{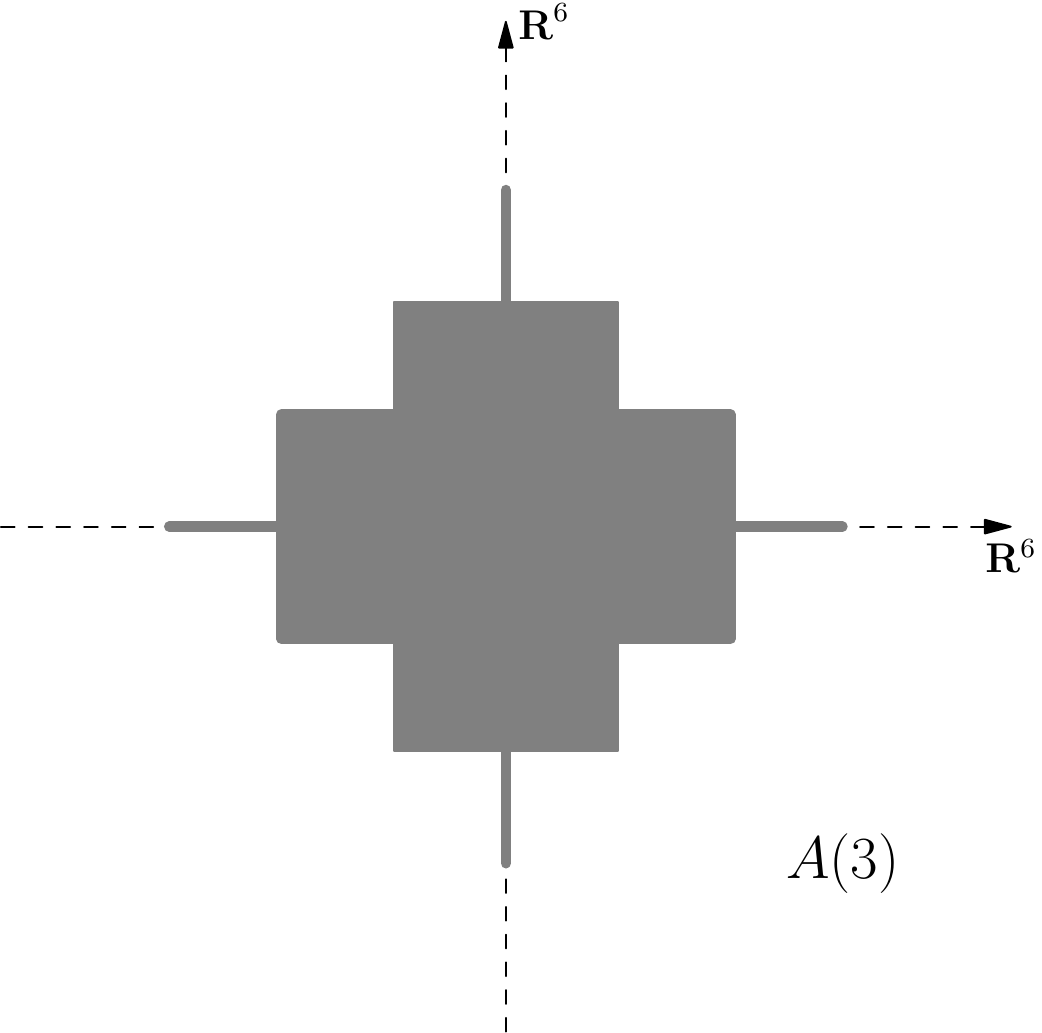}
\end{center}
\caption{A counterexample in $\R^{12}$.}
\end{figure}

For $\R^1$, the conjecture  may be proved by adapting a proof of \cite{GMR10}
on cardinality of integer sumsets; this was also independently observed by F.~Barthe.
Let $k \geq 1$. Set $S=A_1+ \cdots +A_k$ and for  $i \in[k]$, let $a_i=\min A_i$, $b_i=\max A_i$, 
$$S_i = \sum_{j\in[k]\setminus\{i\}}A_j,$$
$s_i=\sum_{j<i}a_j+\sum_{j>i}b_j$, $S_i^-=\{x\in S_i ; x\le s_i\}$ and $S_i^+=\{x\in S_i ; x> s_i\}$. For all $i \in [k-1]$, one has 
$$ S\supset (a_i+S_i^-)\cup(b_{i+1}+S_{i+1}^+).$$
Since $a_i+s_i=\sum_{j\le i}a_j+\sum_{j>i}b_j=b_{i+1}+s_{i+1}$, the above union is a disjoint union. Thus for $i \in [k-1]$
$$\vol_1(S) \ge \vol_1(a_i+S_i^-)+\vol_1(b_{i+1}+S_{i+1}^+) 
= \vol_1(S_i^-) + \vol_1(S_{i+1}^+).$$
Notice that $S_1^-=S_1$ and $S_k^+=S_k\setminus\{s_k\}$, thus adding the above $k-1$ inequalities we obtain
\begin{eqnarray*}
(k-1)\vol_1(S) & \ge & \sum_{i=1}^{k-1}\left(\vol_1(S_i^-) + \vol_1(S_{i+1}^+) \right) 
= \vol_1(S_1^-) + \vol_1(S_k^+) + \sum_{i=2}^{k-1}\vol_1(S_i) \\
&=& \sum_{i=1}^k\vol_1(S_i).
\end{eqnarray*}
Now taking all the sets $A_i=A$, and dividing through by $k(k-1)$, we see that 
we have established Conjecture \ref{weakconj} in dimension 1.

\section{The behavior of Schneider's non-convexity index and the Hausdorff distance}\label{sec:c}

We establish Theorems \ref{thm:c-quant-mono} and \ref{thm:d-quant-mono} in this section. This relies crucially on the elementary observations that $\conv (A+B)=\conv(A)+\conv(B)$
and $(t+s)\conv(A)=t\conv(A) +s\conv(A)$ for any $t,s>0$ and any compact sets $A, B$.\\

{\bf Proof of Theorem \ref{thm:c-quant-mono}.}
Denote $\lambda=c\left(A(k)\right)$. Since $\conv(A(k))=\conv(A)$, from the definition of $c$, one knows that 
$A(k)+\lambda\conv(A)=\conv(A)+\lambda\conv(A)=(1+\lambda)\conv(A)$. Using that $A(k+1)=\frac{A}{k+1}+\frac{k}{k+1}A(k)$, one has 
\begin{eqnarray*}
A(k+1)+\frac{k}{k+1}\lam\conv(A) &=& \frac{A}{k+1}+\frac{k}{k+1}A(k)+\frac{k}{k+1}\lam\conv(A)\\
&=&\frac{A}{k+1}+\frac{k}{k+1}\conv(A)+\frac{k}{k+1}\lam\conv(A)\\
&\supset&\frac{\conv(A)}{k+1}+\frac{k}{k+1}A(k)+\frac{k}{k+1}\lam\conv(A)\\
&=&\frac{\conv(A)}{k+1}+\frac{k}{k+1}(1+\lam)\conv(A)\\
&=&\left(1+\frac{k}{k+1}\lam\right)\conv(A).
\end{eqnarray*}
Since the other inclusion is trivial, we deduce that $A(k+1)+\frac{k}{k+1}\lam\conv(A)$ is convex which proves that 
$$c(A(k+1))\le \frac{k}{k+1}\lam=\frac{k}{k+1}c\left(A(k)\right).$$

{\bf Proof of Theorem \ref{thm:d-quant-mono}.}
Let $k\ge c(A)$, then, from the definitions of $c(A)$ and $d(A(k))$, one has 
\begin{eqnarray*}
\conv(A)=\frac{A}{k+1}+\frac{k}{k+1}\conv(A)&\subset& \frac{A}{k+1}+\frac{k}{k+1}\left(A(k)+d(A(k))B_2^n\right)\\
&=&A(k+1)+\frac{k}{k+1}d(A(k))B_2^n.
\end{eqnarray*}
We conclude that 
$$ d\left(A(k+1)\right)\le\frac{k}{k+1}d\left(A(k)\right). $$

\section{Discussion}

\begin{enumerate}
\item By repeated application of Theorem \ref{thm:c-quant-mono}, it is clear that 
the convergence of $c(A(k))$ is at a rate $O(1/k)$ for any compact set $A \subset \R^n$;
this observation appears to be new. In \cite{FMMZ15}, we study the question of the 
monotonicity of $A(k)$, as well as convergence rates, when considering several
different ways to measure non-convexity, including some not mentioned in this note.

\item Some of the results in this note are of interest when one is considering 
Minkowski sums of different compact sets, not just sums of $A$ with copies of itself.
Indeed, the original conjecture of \cite{BMW11} was of this form, and would
have provided a strengthening of the classical Brunn-Minkowski inequality for
more than 2 sets; of course, that conjecture is false since the weaker Conjecture \ref{weakconj}
is false. Nonetheless we do have some related observations in \cite{FMMZ15}; for instance,
it turns out that in general dimension, for compact sets $A_1, \ldots, A_k$,
\ben
\vol_n\left(\sum_{i=1}^kA_i\right) \ge \frac{1}{k-1}\sum_{i=1}^k\vol_n\left(\sum_{j\in[k]\setminus\{i\}}A_j\right).
\een
For convex sets $B_i$, an even stronger fact is true (that this is stronger may not be immediately obvious,
but if follows from well known results, see, e.g., \cite{MT10}):
\ben
\vol_n(B_1+B_2+B_3) + \vol_n(B_1) \geq \vol_n(B_1+B_2) + \vol_n(B_1+B_3)  .
\een

\item There is a variant of the strong monotonicity of Schneider's index when dealing with different sets.
If $A, B, C$ are subsets of $\R^n$, then it is shown in \cite{FMMZ15} (by a similar argument to that used for Theorem~3)
that
$c(A+B+C) \leq \max\{c(A+B), c(B+C)\}$.

\end{enumerate}

\section*{Acknowledgements}
We are indebted to Fedor Nazarov for valuable discussions, in particular for help in the construction of the counterexample in Theorem \ref{voldef}. 
We also thank Franck Barthe, Dario Cordero-Erausquin, Uri Grupel, Joseph Lehec, Bo'az Klartag and Paul-Marie Samson for interesting discussions.

\end{document}